\nonstopmode
\documentclass[letterpaper,10pt]{amsart}
\usepackage{latexsym}
\usepackage{fancyhdr}
\usepackage{amsmath, amssymb}
\usepackage[ansinew]{inputenc}
\usepackage[all]{xy}
\usepackage{verbatim}
\usepackage{hyperref}
\usepackage{subfigure}
\swapnumbers
\theoremstyle{plain}
\newtheorem*{lemma*}{Lemma}
\newtheorem{lemma}[subsection]{Lemma}
\newtheorem*{theorem*}{Theorem}
\newtheorem{theorem}[subsection]{Theorem}
\newtheorem*{proposition*}{Proposition}

\newtheorem*{corollary*}{Corollary}

\newtheorem*{claim*}{Claim}
\theoremstyle{definition}
\newtheorem*{definition*}{Definition}

\newtheorem*{example*}{Example}
\newtheorem{example}[subsection]{Example}
\newtheorem{examples}[subsection]{Examples}
\newtheorem*{algorithm*}{Algorithm}
\newtheorem*{remark*}{Remark}
\newtheorem{remark}[subsection]{Remark}

\newenvironment{demo}[1]{\par\smallskip\noindent{\bf #1.}}{\par\smallskip}
\numberwithin{equation}{subsection}

\def\al{\alpha}

\def\ga{\gamma}
\def\de{\delta}
\def\ep{\epsilon}

\def\ze{\zeta}

\def\la{\lambda}

\def\si{\sigma}

\def\ta{\tau}

\def\om{\omega}

\def\Ga{\Gamma}
\def\De{\Delta}

\def\La{\Lambda}

\def\Ph{\Phi}

\def\C{\mathbb{C}}

\def\N{\mathbb{N}}

\def\R{\mathbb{R}}

\def\cM{\mathcal{M}}

\def\cR{\mathcal{R}}

\def\p{\partial}

\def\<{\langle}
\def\>{\rangle}
\renewcommand{\o}{\circ}

\let\on=\operatorname

\title[Smooth roots of hyperbolic polynomials]
{Smooth roots of hyperbolic polynomials with definable coefficients}

\author[A. Rainer]
{Armin Rainer}

\address{Armin Rainer: Department of Mathematics, University of Toronto, 
40 St.\ George Street, Toronto, Ontario, Canada M5S 2E4}

\email{armin.rainer@univie.ac.at}

\begin{document}

\begin{abstract} 
We prove that the roots of a definable $C^\infty$ curve of monic hyperbolic polynomials admit a definable $C^\infty$ parameterization,
where `definable' refers to any fixed o-minimal structure on $(\R,+,\cdot)$.
Moreover, we provide sufficient conditions, in terms of the differentiability of the coefficients and the order of contact of the roots, 
for the existence of $C^p$ (for $p \in \N$) arrangements of the roots in both the definable and the non-definable case.
These conditions are sharp in the definable and under an additional assumption also in the non-definable case.
In particular, we obtain a simple proof of Bronshtein's theorem in the definable setting.
We prove that the roots of definable $C^\infty$ curves of complex polynomials can be desingularized by means 
of local power substitutions $t \mapsto \pm t^N$.
For a definable continuous curve of complex polynomials we show that any continuous choice of roots is actually locally absolutely continuous.
\end{abstract}

{\noindent{\small\rm To appear in Israel J.\ Math.}}

\thanks{The author was supported by the Austrian Science Fund (FWF), Grant J2771}
\keywords{hyperbolic polynomials, smooth roots, o-minimality}
\subjclass[2000]{26C10, 30C15, 03C64}
\dedicatory{Dedicated to Peter W.\ Michor on the occasion of his 60th birthday}
\date{April 27, 2009}

\maketitle

\section{Introduction}

A monic polynomial $P(x) = x^n + \sum_{j=1}^n (-1)a_j x^{n-j}$ is called hyperbolic if all its roots are real. 
The study of the regularity of its roots, when $P$ depends smoothly on a real parameter, is a classical topic with important applications 
in PDE and perturbation theory. Rellich \cite{Rellich37I} showed that a real analytic curve of hyperbolic polynomials $P$ admits real analytic roots. 
However, the roots of a $C^\infty$ curve $P$ do in general not allow $C^\infty$ (more precisely, $C^{1,\al}$ for any $\al > 0$) parameterizations. 
All counter-examples (e.g.\ in \cite{Glaeser63R}, \cite{AKLM98}, \cite{BBCP06}) are \emph{oscillating}, meaning that some
derivative switches sign infinitely often near some point, where the multiplicity of the roots changes.
By \cite{AKLM98}, $P$ allows $C^\infty$ roots, if no two roots meet of infinite order.

We show in this note that definability of the coefficients guarantees $C^\infty$ solvability of $C^\infty$ curves of hyperbolic polynomials.
By `definable' we mean definable in some fixed, but arbitrary, o-minimal structure $\cM$ on $(\R,+,\cdot)$.
Definability excludes oscillation, however, infinitely flat functions may be definable in some $\cM$. 
We also provide sufficient conditions, in terms of the differentiability of the coefficients and the order of contact of the roots, 
for the existence of $C^p$ (for $p \in \N$) arrangements of the roots in both the definable and the non-definable case.
These conditions are sharp in the definable and under an additional assumption (automatically satisfied if $n\le 4$) also in the non-definable case.
In particular, we give a simple proof of Bronshtein's theorem in the special case of definable coefficients: 
$C^n$ curves $P$ admit $C^1$ roots 
(see \cite{Bronshtein79}, \cite{Wakabayashi86}, and \cite{ColombiniOrruPernazza08}).
As a consequence $C^{2n}$ curves $P$ admit twice differentiable roots (see \cite{KLM04} and \cite{ColombiniOrruPernazza08}).
Bronshtein's theorem is quite delicate and only poorly understood.

Our results complete the perturbation theory for hyperbolic polynomials. Analogous questions for several parameters require additional assumptions
and are not treated in this paper: The roots of $P(t_1,t_2)(x) = x^2 -(t_1^2+t_2^2)$, for $t_1,t_2 \in \R$, cannot be differentiable at $t_1=t_2=0$. 

If the hyperbolicity assumption is dropped, we cannot hope for parameterizations of the roots satisfying a local Lipschitz condition, 
even if the coefficients are real analytic. 
We prove that the roots of definable $C^\infty$ curves of complex polynomials can be desingularized by means 
of local power substitutions $t \mapsto \pm t^N$.
For definable continuous curves of complex polynomials, we show that any continuous choice of roots is 
actually locally absolutely continuous (not better!). This extends results in \cite{RainerAC}.

I am happy to thank E.\ Bierstone and K.\ Kurdyka for the discussions which led to the writing of this paper. 

\section{Definable functions and smoothness}

\subsection{Multiplicity} \label{defmult}
For a continuous real or complex valued function $f$ defined near $0$ in $\R$, 
let the \emph{multiplicity} $m_0(f)$ at $0$ be the supremum 
of all integers $p$ such that $f(t)=t^p g(t)$ near $0$ for a continuous function $g$. 
Note that, if $f$ is of class $C^n$ and $m_0(f) < n$, then $f(t) = t^{m_0(f)} g(t)$ 
near $0$, where now $g$ is $C^{n-m_0(f)}$ and $g(0) \ne 0$. 
Similarly, one can define the multiplicity $m_t(f)$ of a function $f$ at any $t \in \R$.

\begin{lemma} \label{C^p-ext}
Let $I \subseteq \R$ be an open interval containing $0$. Let $f \in C^0(I,\R)$ and $p \in \N$ such that: 
\begin{enumerate}
\item[(1)] $m_0(f) \ge p$
\item[(2)] $f|_{I \setminus \{0\}} \in C^{p+1}(I \setminus \{0\})$
\item[(3)] $0$ is not an accumulation point of $\p \{t \in I\setminus \{0\} : f^{(p+1)}(t) = 0\}$ 
(where $\p A:=\overline A \setminus  A^\o$ denotes the boundary of $A$) .
\end{enumerate}
Then $f \in C^p(I)$.
\end{lemma}

\begin{demo}{Proof}
We use induction on $p$. Let us assume that the assertion is proved for non-negative integers $< p$. 
Note that (3) implies: 
\begin{enumerate}
\item[(3')] $0$ is not an accumulation point of $\p \{t \in I \setminus \{0\} : f^{(q)}(t) = 0\}$, for any integer $0 \le q \le p+1$. 
\end{enumerate}
So we may suppose that $f \in C^{p-1}(I)$, and, by (1), $f^{(q)}(0)=0$ for $0 \le q \le p-1$. We will show that $f \in C^p(I)$.

Let $t>0$. By (3'), either $f^{(p)} = 0$ identically, or $f^{(p-1)}$ is strictly monotonic for small $t$, say, $t < \de$. 
In the first case $f^{(p)}$ extends continuously to $0$. Consider the second case.
Without loss we may suppose that 
\begin{equation} \label{mono}
f^{(p-1)}(s) > f^{(p-1)}(t) \quad \text{ if  } 0 < s < t< \de 
\end{equation}
(otherwise consider $-f^{(p-1)}$).
Then $f^{(p-1)}(s)/s > f^{(p-1)}(t)/t$ if $0 < s < t < \de$. So 
\[
\lim_{t \searrow 0} \frac{f^{(p-1)}(t)}{t} = \sup_{0<t<\de} \frac{f^{(p-1)}(t)}{t} =: a \in \R \cup \{+\infty\}.
\]
By Taylor's formula, for each $t>0$ there is a $0<\xi(t)<t$ such that 
\[
f(t) = t^{p-1} \cdot \frac{f^{(p-1)}(\xi(t))}{(p-1)!}.
\]
By \eqref{mono}, we have $f^{(p-1)}(\xi(t)) > f^{(p-1)}(t)$, and, thus,
\[
\frac{f^{(p-1)}(t)}{t} < (p-1)! \cdot \frac{f(t)}{t^p}. 
\]
By (1), the right-hand side is convergent as $t \searrow 0$.
So $a < + \infty$. 

By (3), $f^{(p)}$ is strictly monotonic for small $t$, say, $t < \ep$. 
We may conclude that $\lim_{t \searrow 0} f^{(p)}(t)$ is given by either $\sup_{0<t<\ep} f^{(p)}(t)$
or $\inf_{0<t<\ep} f^{(p)}(t)$. 
By Taylor's formula, for each $n \in \N_{>0}$, there is a $0 < \nu(n) < 1/n$ such that
\[
f^{(p)}(\nu(n)) = p! \cdot \frac{f(\tfrac{1}{n})}{(\tfrac{1}{n})^p} = p! \cdot g(\tfrac{1}{n}) \to p! \cdot g(0) \quad \text{ as } n \to \infty,
\]
where $g(t):=f(t)/t^p$ is continuous by (1).
Hence, $\lim_{t \searrow 0} f^{(p)}(t) = p! \cdot g(0)$. By the mean value theorem, we obtain
\[
a = \lim_{n \to \infty} \frac{f^{(p-1)}(\tfrac{1}{n})}{\tfrac{1}{n}} = \lim_{n \to \infty} f^{(p)}(\ze(n)) = p! \cdot g(0),
\]
where $0 < \ze(n) < 1/n$.
(Note that, if $f^{(p)} = 0$ identically, then $g(0)=0$.) 

In a similar way one proves that $\lim_{t \nearrow 0} f^{(p-1)}(t)/t = \lim_{t \nearrow 0} f^{(p)}(t) = p! \cdot g(0)$.
So $f \in C^p(I)$.
\qed\end{demo}

\begin{example}
Note that condition (3) in lemma \ref{C^p-ext} is necessary:
The function $f(t) := e^{-1/t^2} \sin^2 (e^{1/t^4})$,  $f(0):=0$, satisfies $m_0(f)=\infty$ and is $C^\infty$ off $0$, 
but it is not $C^1$ in any neighborhood of $0$.
\end{example}

\subsection{Definable functions} \label{omin}
Cf.\ \cite{vandenDries98}.
Let $\cM=\bigcup_{n \in \N_{>0}} \cM_n$, where each $\cM_n$ is a family of subsets of $\R^n$. 
We say that $\cM$ is an \emph{o-minimal structure} on $(\R,+,\cdot)$ if the following conditions are satisfied:
\begin{enumerate}
\item[(1)] Each $\cM_n$ is closed under finite set-theoretical operations.
\item[(2)] If $A \in \cM_n$ and $B \in \cM_m$, then $A \times B \in \cM_{n+m}$.
\item[(3)] If $A \in \cM_{n+m}$ and $\pi : \R^{n+m} \to \R^n$ is the projection on the first $n$ coordinates, then $\pi(A) \in \cM_m$.
\item[(4)] If $f,g_1,\ldots,g_l \in \R[X_1,\ldots,X_n]$, then $\{x \in \R^n : f(x)=0,g_1(x)>0,\ldots,g_l(x)>0\} \in \cM_n$.
\item[(5)] $\cM_1$ consists of all finite unions of open intervals and points.
\end{enumerate}

For a fixed o-minimal structure $\cM$ on $(\R,+,\cdot)$, we say that $A$ is \emph{$\cM$-definable} if $A \in \cM_n$ for some $n$.
A mapping $f : A \to \R^m$, where $A \subseteq \R^n$, is called \emph{$\cM$-definable} if its graph is $\cM$-definable.

\textbf{From now on let $\cM$ be some fixed, but arbitrary, o-minimal structure on $(\R,+,\cdot)$. If we write \emph{definable} we will always mean $\cM$-definable.} 

\begin{lemma} \label{C^p-def}
Let $I \subseteq \R$ be an open interval containing $0$, let $f : I \to \R$ be definable, and $p,m \in \N$.
\begin{enumerate}
\item If $f \in C^0(I)$ and $m_0(f) \ge p$, then $f$ is $C^p$ near $0$.
\item If $f \in C^p(I)$, then $h(t):=t^m f(t)$ is $C^{p+m}$ near $0$. 
\end{enumerate}
\end{lemma}

\begin{demo}{Proof}
(1) follows from lemma \ref{C^p-ext} and the Monotonicity theorem (e.g.\ \cite{vandenDries98}).

(2) We use induction on $m$. The statement for $m=0$ is trivial. Suppose that $m > 0$.  
By induction hypothesis, $g(t):=t^{m-1} f(t)$ belongs to $C^{p+m-1}(I)$ and
$h^{(p+m-1)}(t) = t g^{(p+m-1)}(t) + (p+m-1) g^{(p+m-2)}(t)$. 
Thus 
\[
\lim_{t \to 0} \frac{h^{(p+m-1)}(t)-h^{(p+m-1)}(0)}{t} = (p+m) g^{(p+m-1)}(0).
\]
Let $t > 0$. By definability, $h^{(p+m)}(t)$ exists and is either a constant $a$ or strictly monotonic for small $t$, say, $t < \ep$. 
Hence, $\lim_{t \searrow 0} h^{(p+m)}(t)$ is given by either $a$, $\sup_{0<t<\ep} h^{(p+m)}(t)$,
or $\inf_{0<t<\ep} h^{(p+m)}(t)$.
By the mean value theorem, for each $n \in \N_{>0}$, there is a $0 < \nu(n) < 1/n$ such that
\[
h^{(p+m)}(\nu(n)) = \frac{h^{(p+m-1)}(\tfrac{1}{n})-h^{(p+m-1)}(0)}{\tfrac{1}{n}} \to (p+m) g^{(p+m-1)}(0) \quad \text{ as } n \to \infty.
\]
So $\lim_{t \searrow 0} h^{(p+m)}(t) = (p+m) g^{(p+m-1)}(0)$. Similarly for $t < 0$.
\qed\end{demo}

\begin{examples}
The conditions in lemma \ref{C^p-def} are sharp:
Let 
\begin{equation} \label{f_p}
f_p(t) := \left\{
\begin{matrix}
t^{p+1} & \text{for } t \ge 0\\ 0 & \text{for } t <0
\end{matrix}
\right..
\end{equation}
Then $m_0(f_p) = p$, and $f_p$ is $C^{p,1}$ but not $C^{p+1}$.
Moreover, $f_{p+m}(t) = t^m f_p(t)$ is $C^{p+m,1}$ but not $C^{p+m+1}$.
\end{examples}

\section{Smooth square roots}

\subsection{}

Let $I \subseteq \R$ be an open interval.
If $f : I \to \R_{\ge 0}$ is definable and continuous, then
$\{t \in I : 0 < m_t(f) < \infty\} \subseteq \p \{t \in I : f(t)=0\}$.
So
\[
2 \overline m(f) := \sup \{m_t(f) < \infty : t \in I\}
\]
is a well-defined integer. 
If $f$ is $C^n$ and $n > 2 \overline m(f)$, then $\overline m(f)$ is the maximal \emph{finite} order of vanishing of the square roots of $f$.

\begin{theorem} \label{n=2}
Let $I \subseteq \R$ be an open interval, $f: I \to \R_{\ge 0}$ a non-negative definable function, and $p \in \N_{>0}$. 
Consider $P(t)(x) = x^2 -f(t)$. 
Then we have:
\begin{enumerate}
\item[(1)] If $f$ is $C^\infty$, then the roots of $P$ admit definable $C^\infty$ parameterizations.
\item[(2)] If $f$ is $C^{p+2 \overline m(f)}$, then the roots of $P$ admit definable $C^{p+\overline m(f)}$ parameterizations.
\end{enumerate}
\end{theorem}

\begin{demo}{Proof}
We prove (1) and (2) simultaneously and indicate differences when arising.

Note that any continuous choice of roots is definable (cf.\ lemma \ref{defroots}).

Let $t_0 \in I$. If $0 \le m_{t_0}(f)< \infty$, then $m_{t_0}(f) = 2m$ for some $m \in \N$, since $p+2 \overline m(f)-m_{t_0}(f)\ge 1$ and $f \ge 0$. 
So $f(t) = (t-t_0)^{2m} f_{(m)}(t)$, where
\[
f_{(m)}(t) = \int_0^1 \frac{(1-r)^{2m-1}}{(2m-1)!} f^{(2m)}(t_0 + r (t-t_0)) dr
\]
is $C^\infty$ (resp.\ $C^{p+2 \overline m(f) - 2 m}$), definable, and $f_{(m)}(t_0) >0$. . 
Then the functions $g_{\pm}(t) : = \pm (t-t_0)^m \sqrt{f_{(m)}(t)}$ are $C^\infty$ (resp.\ $C^{p+2 \overline m(f) - m}$, by lemma \ref{C^p-def}(2))
and represent the roots of $P$ near $t_0$.

Now assume that $m_{t_0}(f) = \infty$. 
In a neighborhood of $t_0$, consider the continuous functions $g_\pm(t):=\pm \sqrt{f(t)}$. 
Then $m_{t_0}(g_\pm) = \infty$.  
By lemma \ref{C^p-def}(1), for each $p$, there is a neighborhood $I_p$ of $t_0$ such that the roots $g_\pm$ are $C^p$ on $I_p$. 
Now, either $f=0$ identically near $t_0$, 
or $t_0$ belongs to $\p (f^{-1}(0))$ which is finite, by definability. 
Thus, in case (1), $g_\pm$ is $C^\infty$ off $t_0$, and, hence, near $t_0$.

So for each $t_0 \in I$ we have found local $C^\infty$ (resp.\ $C^{p+\overline m(f)}$) parameterizations of the roots of $P$ near $t_0$. 
One can glue these to a global parameterization, see \ref{main}(4) below.
\qed\end{demo}

\begin{examples}
The condition in theorem \ref{n=2}(2) is sharp:
The non-negative function $f(t) := t^{2m} (1 + f_p(t))$, where $f_p$ is defined in \eqref{f_p}, is $C^{p+2m,1}$ but not $C^{p+2m+1}$.
Its square roots $g_{\pm}(t) : = \pm t^m \sqrt{1 + f_p(t)}$ are $C^{p+m}$ but not $C^{p+m+1}$.
\end{examples}

\section{Smooth roots of hyperbolic polynomials}

\subsection{} \label{pre} 
Let
\[
P(z) = z^n + \sum_{j=1}^n (-1)^j a_j z^{n-j} = \prod_{j=1}^n (z-\la_j)
\]
be a monic polynomial with complex coefficients $a_1,\ldots,a_n$ 
and roots $\la_1,\ldots,\la_n$. 
By Vieta's formulas, $a_i = \sigma_i(\la_1,\ldots,\la_n)$, where $\sigma_1,\ldots,\sigma_n$ 
are the elementary symmetric functions in $n$ variables:
\begin{equation} \label{esf}
\si_i(\la_1,\ldots,\la_n) = 
\sum_{1 \le j_1 < \cdots < j_i \le n} \la_{j_1} \cdots \la_{j_i}.
\end{equation}
Denote by $s_i$, $i \in \N$, the Newton polynomials 
$\sum_{j=1}^n \la_j^i$ which are related to the elementary symmetric functions 
by
\begin{equation} \label{rec}
s_k - s_{k-1} \sigma_1 + s_{k-2} \sigma_2 - \cdots 
+ (-1)^{k-1} s_1 \sigma_{k-1} + (-1)^k k \sigma_k = 0, \quad (k \ge 1).
\end{equation}

Let us consider the so-called Bezoutiant
\[ 
B := 
\begin{pmatrix} 
s_0 & s_1 & \ldots & s_{n-1}\\ 
s_1 & s_2 & \ldots & s_n \\ 
\vdots & \vdots & \ddots & \vdots\\ 
s_{n-1} & s_n & \ldots &  s_{2n-2} 
\end{pmatrix} 
= \left(s_{i+j-2}\right)_{1 \le i,j \le n}.
\]
Since the entries of $B$ are symmetric polynomials in $\la_1,\ldots,\la_n$, 
we find a unique symmetric $n \times n$ matrix $\tilde{B}$ with 
$B = \tilde{B} \circ \sigma$, where $\sigma =(\sigma_1,\ldots,\sigma_n)$. 

Let $B_k$ denote the minor formed by the first $k$ rows and columns of $B$. Then we have 
\begin{equation} \label{eqdel}
\Delta_k(\la) := \on{det} B_k(\la) = \sum_{i_1 < \cdots < i_k} 
(\la_{i_1}-\la_{i_2})^2 \cdots (\la_{i_1}-\la_{i_k})^2 \cdots (\la_{i_{k-1}}-\la_{i_k})^2.
\end{equation}
Since the polynomials $\Delta_k$ are symmetric, we have 
$\Delta_k = \tilde{\Delta}_k \circ \si$ 
for unique polynomials $\tilde{\Delta}_k$.

By \eqref{eqdel}, the number of
distinct roots of $P$ equals the
maximal $k$ such that $\tilde \De_k(P) \ne 0$.
(Abusing notation we identify $P$ with the $n$-tupel $(a_1,\ldots,a_n)$ of its coefficients when convenient.)

If all roots $\la_j$ (and thus all coefficients $a_j$) of $P$ are real, we say that $P$ is \emph{hyperbolic}.

\begin{theorem*}[Sylvester's version of Sturm's theorem, see e.g.\ \cite{Procesi78} for a modern proof]
Suppose that all coefficients of $P$ are real.
Then $P$ is hyperbolic if and only if $\tilde{B}(P)$ is
positive semidefinite. The rank of $\tilde{B}(P)$ 
equals the number of distinct roots
of $P$ and its signature equals the number of distinct real roots.
\end{theorem*}

\begin{lemma}[Splitting lemma {\cite[3.4]{AKLM98}}]  \label{split}
Let $P_0 = z^n + \sum_{j=1}^n (-1)^j a_j z^{n-j}$ be a polynomial satisfying 
$P_0 = P_1 \cdot P_2$, where $P_1$ and $P_2$ are polynomials without common root. 
Then for $P$ near $P_0$ we have $P = P_1(P) \cdot P_2(P)$ 
for analytic mappings 
of monic polynomials $P \mapsto P_1(P)$ and $P \mapsto P_2(P)$, 
defined for $P$ 
near $P_0$, with the given initial values.
\end{lemma}

\subsection{} \label{set}
For the rest of the section,
let $I \subseteq \R$ be an open interval
and consider a (continuous) curve of hyperbolic polynomials
\[
P(t)(x) = x^n + \sum_{j=1}^n (-1)^j a_j(t) x^{n-j}, \quad (t \in I).
\]
Then the roots of $P$ admit a continuous parameterization, e.g., ordering them by size,
$\la_1 \le \la_2 \le \cdots \le \la_n$.

\begin{lemma} \label{defroots}
If the coefficients $a_j$ of $P$ are definable, then every continuous parameterization $\la_j$ of the roots of $P$ is definable. 
\end{lemma}

\begin{demo}{Proof}
Ordering the roots of $P$ by size, $\mu_1 \le \mu_2 \le \cdots \le \mu_n$, gives a continuous parameterization which is 
evidently definable. Since all $\tilde \De_k \o P$ are definable, 
the set $E$ of $t \in I$ where the multiplicity of the roots changes is finite. 
The complement of $E$ consists of finitely many intervals, on each of which
the parameterizations $\la_j$ and $\mu_j$ differ only by a constant permutation. 
Thus each $\la_j$ is definable. 
\qed\end{demo}

\begin{lemma}[Multiplicity lemma {\cite[3.7]{AKLM98}}]   \label{mult}
Suppose that $0 \in I$ and that $a_1=0$ identically.
Let $r \in \N$. 
If each $a_j \in C^{nr}(I)$, 
then the following conditions are equivalent:
\begin{enumerate}
\item[(1)] $m_0(a_k) \ge k r$, for all $2 \le k \le n$.
\item[(2)] $m_0(\tilde{\Delta}_k) \ge k (k-1) r$, for all $2 \le k \le n$.
\item[(3)] $m_0(a_2) \ge 2 r$.
\end{enumerate}
\end{lemma}

\begin{demo}{Proof}
Obvious modification of the proof of \cite[3.7]{AKLM98}.
\qed\end{demo}

\subsection{} \label{Einfty}
Let $E^{(\infty)}(P)$ denote the set of all $t \in I$ which satisfy following condition:
\begin{enumerate}
\item[(\#)] Let $s=s(t,P)$ be maximal with the property that the germ at $t$ of $\tilde \De_s \o P$ is not $0$. 
Then $m_t(\tilde \De_s \o P) = \infty$.
\end{enumerate}
Consider the condition:
\begin{enumerate}
\item[(\#')] There exists a continuous parameterization $\la_j$ of the roots of $P$ 
such that distinct $\la_j$ meet of infinite order at $t$, i.e.,
there exist $i \ne j$ such that the germs of $\la_i$ and $\la_j$ at $t$ do not coincide and $m_t(\la_i-\la_j)=\infty$.
\end{enumerate}
By \eqref{eqdel}, (\#') implies (\#).

If the coefficients of $P$ (and thus the $\tilde \De_k \o P$) are definable, 
then $E^{(\infty)}(P)$ is finite and the family of continuous parameterizations of the roots of $P$ is finite.
Then (\#) and (\#') are equivalent.

\subsection{} \label{ord}
Let $t_0 \in I$.
Choose a continuous parameterization $\la_j$ of the roots of $P$.
We denote by $\overline m_{t_0}(P,\la)$ the maximal \emph{finite} order of contact of the $\la_j$ at $t_0$, i.e.,
\[
\overline m_{t_0}(P,\la) = \max \{m_{t_0}(\la_i-\la_j) < \infty : 1 \le i < j \le n\}.
\]
The integer $\overline m_{t_0}(P,\la)$ depends on the choice of the $\la_j$. 

If $t_0 \not\in E^{(\infty)}(P)$ and $s = s(t_0,P)$ is the integer defined in (\#), then, by \eqref{eqdel},
\[
\overline m_{t_0}(P,\la) \le \frac{m_{t_0}(\tilde \De_s \o P)}{2}.
\]

If the coefficients of $P$ are definable, 
then the family of continuous parameterizations of the roots of $P$ is finite.

Hence, 
\[
\overline m_{t_0}(P) := \sup_{\la} \overline m_{t_0}(P,\la),
\]
where $\la$ is any continuous arrangement of the roots of $P$, is a well-defined integer,
if either $t_0 \not\in E^{(\infty)}(P)$ or the coefficients of $P$ are definable. 
It is the \emph{maximal finite order of contact of the roots of $P$}.

\begin{lemma} \label{prop}
Suppose that either $t_0 \not\in E^{(\infty)}(P)$ or the coefficients of $P$ are definable.
We have:
\begin{enumerate}
\item[(1)] If $P=P_1 \cdot P_2$ as provided by the splitting lemma \ref{split}, then 
$\overline m_{t_0}(P)=\max \{\overline m_{t_0}(P_1),\overline m_{t_0}(P_2)\}$.
\end{enumerate}
Assume that all roots of $P(t_0)$ coincide. Then:
\begin{enumerate}
\item[(2)] Replacing the variable $x$ with $x-a_1(t)/n$, leaves $\overline m_{t_0}(P)$ unchanged.
\item[(3)] If $a_1=0$, then $m_{t_0}(a_2) \le 2 \overline m_{t_0}(P)+1$.
\item[(4)] Suppose that $a_1=0$ and $a_k(t) =(t-t_0)^{kr} a_{(r),k}(t)$ for continuous $a_{(r),k}$, $2 \le k \le n$, and some $r \in \N_{>0}$.
Consider
\[
P_{(r)}(t)(x) := 
x^n + \sum_{j=2}^n (-1)^j a_{(r),k}(t) x^{n-j}.
\]
Then $\overline m_{t_0}(P_{(r)}) \le \overline m_{t_0}(P) -r$.
\end{enumerate}
\end{lemma}

\begin{demo}{Proof}
(1) and (2) are immediate from the definition.
(3) is a consequence of $-2n a_2 = \sum_{i<j}(\la_i - \la_j)^2$ and the fact that, for a continuous function $f$, 
we have $m_{t_0}(f^2) \le 2 m_{t_0}(f) +1$.  
(4) follows from the observation that, if $t \mapsto \la_i(t)$ parameterize the roots of $t \mapsto P_{(r)}(t)$, then $t \mapsto (t-t_0)^r \la_i(t)$
represent the roots of $t \mapsto P(t)$.
\qed\end{demo}

\begin{example}
Note that in \ref{prop}(3) equality can occur:
Let $f(t) := t^{3+1/3}$ for $t \ge 0$ and $f(t) := 0$ for $t <0$, and consider $P(t)(x) = x^2-f(t)$.
Then $m_0(f) = 3$ and $\overline m_0(P)=1$.
\end{example}

\subsection{} \label{mJ}
If the coefficients $a_j$ of $P$ (and thus the $\tilde \De_k \o P$) are definable, 
then the set $\{t \in I : \overline m_t(P)>0\}$ is finite and
\[
\overline m(P) = \overline m_I(P) := \sup \{\overline m_t(P) : t \in I\}
\]
is a well-defined integer. 

\begin{lemma} \label{tree}
For $n \in \N_{>0}$ let $\mathcal R(n)$ denote the family of all rooted trees $T$ with vertices labeled in the following way: 
the root is labeled $n$, the labels of the successors of a vertex labeled $m$ form a partition of $m$, 
the leaves (vertices with no successors) are all labeled $1$. 
Define $d(n):= \max_{T \in \mathcal R(n)} \{\text{sum over all labels $\ge 2$ in $T$}\}$.
Then 
\begin{equation} \label{dn}
d(n) = \frac{1}{2} n (n+1)-1.
\end{equation}
\end{lemma}

\begin{demo}{Proof}
Observe that $d(1)=0$. 
Then \eqref{dn} is equivalent to $d(n+1)=n+1+d(n)$ for $n \ge 1$.
We use induction on $n$.
It suffices to show $d(n) \ge d(n_1)+\cdots+d(n_p)$ 
for $n_1+\cdots+n_p=n+1$, where $p \ge 2$ and $n_i \in \N_{>0}$. 
By induction hypothesis, this inequality is equivalent to 
\begin{align*}
&\frac{1}{2} n (n+1)-1 \ge \frac{1}{2} n_1 (n_1+1) + \cdots + \frac{1}{2} n_p (n_p+1)-p\\
&\Longleftrightarrow \quad \frac{1}{2} ((n_1+\cdots+n_p)^2-(n_1^2+\cdots+n_p^2)) \ge n_1+\cdots+n_p - p +1\\
&\Longleftarrow \quad (n_1+\cdots+n_{p-1})n_p \ge n_1 + \cdots + n_p - 1.
\end{align*}
The last inequality has the form $ab \ge a+b-1$ for $a,b \in \N_{>0}$, which is easily verified.
\qed\end{demo}

Note that $d(n)+n$ computes the maximal sum of all degrees occurring in a repeated splitting of a polynomial of degree $n$ into 
a product of polynomials of strictly smaller degree until each factor has degree one. 

\begin{theorem} \label{main}
Let $I \subseteq \R$ be an open interval. Consider a curve of hyperbolic polynomials
\[
P(t)(x) = x^n + \sum_{j=1}^n (-1)^j a_j(t) x^{n-j}, \quad (n \ge 2),
\]
with definable coefficients $a_j$.
Let $p \in \N_{>0}$ and $d(n)=n(n+1)/2-1$.
Then:
\begin{enumerate}
\item[(1)] If the $a_j$ are $C^\infty$,
then the roots of $P$ can be parameterized by definable $C^\infty$ functions, globally.
\item[(2)] If the $a_j$ are $C^{p+1+d(n) \overline m(P)}$,
then the roots of $P$ can be parameterized by definable $C^p$ functions, globally.
\end{enumerate}
\end{theorem}

The condition in \ref{main}(2) is not best possible. However, it is convenient to prove this preliminary result parallel to the $C^\infty$ case and 
strengthen it in theorem \ref{sharp} below.

\begin{demo}{Proof}
We prove (1) and (2) simultaneously and indicate differences when arising.
Any continuous parameterization of the roots of $P$ is definable, by lemma \ref{defroots}.

We proceed by induction on $n$. The case $n=2$ is covered by theorem \ref{n=2} (since we may always assume $a_1=0$, see (II) below). 
Suppose the assertion is proved for degrees $<n$.

\begin{claim*}[3]
There exists a local $C^\infty$ (resp.\ $C^p$) parameterization $\la_i$ of the roots of $P$ near each 
$t_0 \in I$. 
The local $C^\infty$ choices $\la_i$ of the roots are unique in the following sense: 
\begin{enumerate}
\item[($\star$)] On the set $\{\la_1,\ldots,\la_n\}$ consider the equivalence relation $\la_i \sim \la_j$ iff $m_{t_0}(\la_i-\la_j)=\infty$.
If $\mu_i$ is a different local $C^\infty$ parameterization of the roots of $P$ near $t_0$, then 
$\{\la_1,\ldots,\la_n\}/\!\!\sim~= \{\mu_1,\ldots,\mu_n\}/\!\!\sim$.
\end{enumerate}
\end{claim*}

Note that ($\star$) is trivially satisfied if $n=2$. 
Without loss we may assume that $0 \in I$ and $t_0=0$. 
We distinguish different cases:

(I) If there are distinct roots at $0$, we may factor $P(t) = P_1(t) \cdot P_2(t)$ in an open subinterval $I_0 \ni 0$ 
such that $P_1$ and $P_2$ have no common roots, by the splitting lemma \ref{split}.
The coefficients of each $P_i$ are definable, since its roots are.
By lemma \ref{prop}(1), we have 
 \[
\overline m_{I_0}(P) = \max \{\overline m_{I_0}(P_1),\overline m_{I_0}(P_2)\}.
\]
By the induction hypothesis, $P_1$ and $P_2$ (and hence $P$) admit $C^\infty$ (resp.\ $C^p$) parameterizations of 
its roots on $I_0$ which are unique in the sense of ($\star$) in case (1).

(II) If all roots of $P(0)$ coincide, then we first reduce $P$ to the case 
$a_1 = 0$, by replacing $x$ by $x - a_1(t)/n$ (which leaves  
$\overline m(P)$ and ($\star$) unchanged, 
by lemma \ref{prop}(2)). 
Then all roots of $P(0)$ are equal to $0$. 
So $a_2(0)=0$.
Clearly, the new coefficients are still definable.

(IIa) If $m_0(a_2)$ is finite, then $p+1+d(n) \overline m(P)-m_0(a_2) \ge 1$, by lemma \ref{prop}(3).
So $m_0(a_2) = 2 r$ for some $r \in \N_{>0}$, since
$0 \le \tilde \De_2 = -2 n a_2$. 
Let 
\[
q:=p+1+d(n) \overline m(P)-nr.
\]  
By the multiplicity lemma \ref{mult}, we obtain $a_k(t) = t^{k r} a_{(r),k}(t)$ 
for definable $C^\infty$ 
(resp.\ $C^q$) functions $a_{(r),k}$ and $2 \le k \le n$. 
Consider the $C^\infty$ (resp.\ $C^q$) curve of 
hyperbolic polynomials 
\begin{equation} \label{Pr}
P_{(r)}(t)(x) := 
x^n + \sum_{j=2}^n (-1)^j a_{(r),k}(t) x^{n-j}.
\end{equation}
Since $a_{(r),2}(0) \ne 0$, not all roots of $P_{(r)}(0)$ coincide. 
We have $d(n)-n=d(n-1)$ and, by lemma \ref{prop}(4), $\overline m(P_{(r)}) \le \overline m(P)-r$.
Thus, the splitting lemma \ref{split} and the induction hypothesis
provide $C^\infty$ (resp.\ $C^p$) parameterizations $\la_j$ of the roots of $P_{(r)}$ near $0$ which are unique in the sense of ($\star$). 
But then the $C^\infty$ (resp.\ $C^p$) functions $t \mapsto t^r \la_j(t)$ represent the roots of $t \mapsto P(t)$ near $0$ 
and they are unique in the sense of ($\star$) in case (1).

(IIb) If $m_0(a_2)=\infty$ and $a_2 = 0$, then all roots of $P$ are identically $0$.

(IIc) Finally, if $m_0(a_2)=\infty$ and $a_2 \ne 0$, then, since $-2a_2=\sum_{j=1}^n \la_j^2$,
for any continuous choice of the roots $\la_j$ we find $m_0(\la_j)=\infty$ (for all $j$).
By lemma \ref{C^p-def}(1), for each $p$, there is a neighborhood $I_p$ of $0$ such that the roots $\la_j$ are $C^p$ on $I_p$.
Since $a_2$ is definable, 
for small $t \ne 0$ either not all $\la_j(t)$ coincide or all $\la_j$ are identically $0$ (to the left or the right of $0$).
So, in case (1), all $\la_j$ are $C^\infty$ off $0$, by the splitting lemma \ref{split} and the induction hypothesis, and hence also near $0$.

\begin{claim*}[4]
We may glue the local $C^\infty$ (resp.\ $C^p$) parameterizations of the roots to form a global parameterization.
\end{claim*}

In case (1) the local $C^\infty$ choices of the roots of $P$ can be glued by their uniqueness in the sense of ($\star$).
If $C^\infty$ roots meet of infinite order at $t_0$, any permutation on one side of $t_0$ preserves smoothness.

For (2): Let $\la=(\la_1,\ldots,\la_n)$ be a $C^p$ parameterization of the roots of $P$
defined on a maximal open interval $I_1 \subseteq I$. For contradiction, assume that the right (say) endpoint $t_1$ of $I_1$
belongs to $I$. By claim (3), there exists a local $C^p$ parameterization $\mu=(\mu_1\ldots,\mu_n)$ 
of the roots of $P$ near $t_1$.
Let $t_0$ be in the common domain of $\la$ and $\mu$. Consider a sequence $t_k \searrow t_0$.
For each $k$, there is a permutation $\ta_k$ of $\{1,\ldots,n\}$ such that 
$\la(t_k) = \ta_k.\mu(t_k)$.
By passing to a subsequence, we can assume that 
$\la(t_k) = \ta.\mu(t_k)$ for all $k$ and a fixed permutation $\ta$.
Thus, $\la(t) = \ta.\mu(t)$ for all $t\ge t_0$, by definability. 
So $\tilde \la(t) := \la(t)$ for $t \le t_0$ and $\tilde \la(t):=\ta.\mu(t)$ for $t \ge t_0$ defines a 
$C^p$ parameterization on a larger interval, a contradiction.
\qed\end{demo}

\begin{remark}
Suppose that $\overline m(P)=0$. Then the roots of $P$ do not meet or they meet \emph{slowly}, i.e., 
$(\la_i(t)-\la_j(t))/t$ is not continuous at $t=0$. In the latter case $a_2 \not\in C^2$, by \ref{prop}(3),
and so \ref{main}(2) is empty.
\end{remark}

\section{Sharp sufficient conditions for \texorpdfstring{$C^p$}{} roots}

The conditions in theorem \ref{main}(2) are not sharp. 
We shall obtain sharp sufficient conditions for $C^p$ roots, given that the coefficients are definable.
In the non-definable case we still get sharp sufficient conditions, if $P$ is of a special type.
The proof of \ref{main}(2) was not for nothing, since it is needed in the definition of $\Ga$ and $\ga$ below.

\subsection{The definable case} \label{defGaga}
Let $P(t)$, $t \in I$, be a curve of monic hyperbolic polynomials of degree $n$ with definable $C^{d(n) \overline m(P) +2}$ coefficients $a_j$.
For each $t_0 \in I$, let us define two integers $\Ga_{t_0}(P)$ and $\ga_{t_0}(P)$ inductively:

(I) If $P(t) = P_1(t) \cdot P_2(t)$ near $t_0$, and $P_i(t_0)$, $i=1,2$, have distinct roots,  
\begin{align}
\Ga_{t_0}(P) &:= \max\{\Ga_{t_0}(P_1),\Ga_{t_0}(P_2)\}, \label{Ga_split}\\
\ga_{t_0}(P) &:= \Ga_{t_0}(P) - \max\{\Ga_{t_0}(P_1)-\ga_{t_0}(P_1),\Ga_{t_0}(P_2)-\ga_{t_0}(P_2)\}. \label{ga_split}
\end{align}

(II) If $\deg(P) >1$ and all roots of $P(t_0)$ coincide, reduce to the case $a_1=0$ (without changing $\Ga_{t_0}(P)$ and $\ga_{t_0}(P)$). 
If $m_{t_0}(a_2) = 2 r<\infty$, consider $P_{(r)}$ as in \eqref{Pr} (for $t_0$ instead of $0$), and set
\begin{align}
\Ga_{t_0}(P) &:= \Ga_{t_0}(P_{(r)}) + \deg(P) r, \label{Ga_(r)}\\
\ga_{t_0}(P) &:= \ga_{t_0}(P_{(r)})+ r. \label{ga_(r)}
\end{align}
If $m_{t_0}(a_2) = \infty$, set $\Ga_{t_0}(P) := 0$ and $\ga_{t_0}(P) := 0$.

(III) If $\deg(P) = 1$, set $\Ga_{t_0}(P) := 0$ and $\ga_{t_0}(P) := 0$.

Note that, (by the proof of \ref{main}(2)) the coefficients of $P$ being in $C^{d(n) \overline m(P) +2}$, 
guarantees that $\Ga_{t_0}(P)$ and $\ga_{t_0}(P)$ are well-defined.
With hindsight it suffices to assume that the coefficients of $P$ belong to $C^{\Ga_{t_0}(P) + 1}$ near $t_0$.

Since the coefficients of $P$ are definable, the set of $t_0 \in I$ such that $\Ga_{t_0}(P)>0$ or $\ga_{t_0}(P)>0$ is finite and
\begin{align}
\Ga(P) &:= \sup \{\Ga_{t_0}(P) : t_0 \in I\}, \label{Ga_sup}\\
\ga(P) &:= \Ga(P) - \sup \{\Ga_{t_0}(P) - \ga_{t_0}(P) : t_0 \in I\} \label{ga_sup}
\end{align}
are well-defined integers. By construction, we have
\[
\ga(P) \le \Ga(P) \le d(n) \overline m(P) +1.
\]

If $P(t)(x)=x^2-f(t)$ (where $f \ge 0$), then $\Ga(P)=2 \overline m(f)$ and $\ga(P) = \overline m(f)$.

\begin{theorem} \label{sharp}
Let $I \subseteq \R$ be an open interval. Consider a curve of hyperbolic polynomials
\[
P(t)(x) = x^n + \sum_{j=1}^n (-1)^j a_j(t) x^{n-j}, 
\]
with definable coefficients $a_j$.
For each $p \in \N_{>0}$, we have:
\begin{enumerate}
\item[(1)] If the $a_j$ are $C^{p+\Ga(P)}$,
then the roots of $P$ can be parameterized by definable $C^{p+\ga(P)}$ functions, globally.
\end{enumerate}
\end{theorem}

\begin{demo}{Proof}
By \ref{main}(4) it suffices to show the local assertion.
Let $t_0 \in I$ be fixed. 

\begin{claim*}[2]
If the $a_j$ are $C^{p+\Ga_{t_0}(P)}$, then the roots of $P$ can be chosen in $C^{p+\ga_{t_0}(P)}$, locally near $t_0$.
\end{claim*}

We use induction on $n$ and follow the steps in \ref{main}. The case $n=1$ is trivial and $n=2$ is treated in theorem \ref{n=2}(2).
Without loss assume that $0 \in I$ and $t_0=0$.

(I) If $P(0)$ has distinct roots, we have a factorization $P(t)=P_1(t) \cdot P_2(t)$ near $0$, by the splitting lemma \ref{split}.
The coefficients of each factor $P_i$ belong to $C^{p+\Ga_{0}(P)}$. Let $p_i := p+\Ga_{0}(P)-\Ga_{0}(P_i)$. Then $p_i \ge p$, by \eqref{Ga_split}.
By the induction hypothesis, the roots of $P_i$ admit a local parameterization in $C^{p_i+\ga_{0}(P_i)}$. 
By \eqref{ga_split}, we obtain $p_i+\ga_{0}(P_i) \ge p+\ga_{0}(P)$, hence claim (2).

(II) If all roots of $P(0)$ coincide, we reduce to the case $a_1=0$. So $a_2(0)=0$.

(IIa) If $m_{0}(a_2) =2r < \infty$, consider $P_{(r)}$ as in \eqref{Pr}.
The coefficients of $P_{(r)}$ are in $C^{p+\Ga_{0}(P)-nr}$ and $a_{(r),2}(0)\ne 0$. By \eqref{Ga_(r)} and (I), 
there are $C^{p+\ga_{0}(P_{(r)})}$ functions $\la_j$ which represent the roots of $P_{(r)}$ near $0$. 
Then $t \mapsto t^r \la_j(t)$ form a local parameterization of the roots of $P$ which is $C^{p+\ga_{0}(P)}$, 
by lemma \ref{C^p-def}(2) and \eqref{ga_(r)}.

(IIb/c) If $m_{0}(a_2)=\infty$, then $m_{0}(\la_j)=\infty$ for each continuous choice of roots $\la_j$, and we are done, by lemma \ref{C^p-def}(1).

\begin{claim*}[3]
If the $a_j$ are $C^{p+\Ga(P)}$, then the roots of $P$ can be chosen in $C^{p+\ga(P)}$, locally near $t_0$.
\end{claim*}

By claim (2), the roots of $P$ can be chosen in $C^{p+\Ga(P)-\Ga_{t_0}(P)+\ga_{t_0}(P)}$, locally near $t_0$.
By \eqref{ga_sup}, we have $p+\Ga(P)-\Ga_{t_0}(P)+\ga_{t_0}(P) \ge p+\ga(P)$.
\qed\end{demo}

\begin{examples}
The condition in theorem \ref{sharp} is sharp:
Let $p \in \N_{> 3}$ and let $f_p$ be the function defined in \eqref{f_p}.
Consider the $C^{p,1}$ curve of polynomials 
\[
P_p(t)(x) = x^3 -f_p(t) x^2 + (2 f_p(t)-t^2) x -f_p(t).
\]
For the discriminant of $P_p$ we find $\tilde \De_3(P_p(t)) = t^6 (4+o(1))$ if $t\ge 0$ (as $p \ge 3$)
and $\tilde \De_3(P_p(t)) = 4 t^6$ if $t < 0$.
Thus, for small $t$, $P_p(t)$ is hyperbolic.
It is easy to compute $\Ga(P_p)=3$($<p$) and $\ga(P_p)=1$.
By theorem \ref{sharp}, $P_p$ admits $C^{p-2}$ roots.
Suppose, for contradiction, that $P_p$ has $C^{p-1}$ roots $\la_j$.
Since $m_0(\la_j) \ge 1$, we have $\la_j(t) = t \mu_j(t)$ for $C^{p-2}$ functions $\mu_j$.
But then $f_p(t) = t^3 \mu_1(t) \mu_2(t) \mu_3(t)$ is $C^{p+1}$, by lemma \ref{C^p-def}(2), a contradiction.
\end{examples}

\subsection{The non-definable case}
Let $P(t)$, $t \in I$, be a curve of monic hyperbolic polynomials of degree $n$ (not necessarily definable). Assume $E^{(\infty)}(P)=\emptyset$.
We will prove analogs of theorem \ref{main}(2) and, if $P$ is of a special type, of  theorem \ref{sharp}. 
Without the assumption $E^{(\infty)}(P)=\emptyset$, we cannot hope for $C^{1,\al}$ roots (for any $\al > 0$), even if the coefficients are $C^\infty$
(e.g.\ \cite{Glaeser63R}, \cite{AKLM98}, \cite{BBCP06}).

Let $J \subseteq I$ be a compact subinterval of $I$.
Define
\[
\overline m_J(P) := \sup \{\overline m_t(P) : t \in J\} \in \N \cup \{+\infty\}.
\]
The interesting case is $\overline m_J(P) < \infty$, but what follows is also true for $\overline m_J(P) = \infty$.

Assume that $P$ has $C^{d(n) \overline m_J(P) +2}$ coefficients $a_j$.
For each $t_0 \in I$, we can define the two integers $\Ga_{t_0}(P)$ and $\ga_{t_0}(P)$ in the same way as in \ref{defGaga}.
Again it is enough to assume that the $a_j$ belong to $C^{\Ga_{t_0}(P) + 1}$ near $t_0$.
Define $\Ga_J(P), \ga_J(P) \in \N \cup \{+\infty\}$ by
\begin{align}
\Ga_J(P) &:= \sup \{\Ga_{t_0}(P) : t_0 \in J\}, \label{Ga_J}\\
\ga_J(P) &:= \Ga_J(P) - \sup \{\Ga_{t_0}(P) - \ga_{t_0}(P) : t_0 \in J\}. \label{ga_J}
\end{align}
By construction, 
\[
\ga_J(P) \le \Ga_J(P) \le d(n) \overline m_J(P) +1.
\]
The interesting case is when $\Ga_J(P)$ and $\ga_J(P)$ are finite, but what follows is true in any case.

\begin{theorem} \label{Cp}
If the coefficients $a_j$ of $P$ are $C^{p + \Ga_J(P)}$,
then the roots of $P$ can be parameterized by $C^p$ functions, globally near $J$.
\end{theorem}

\begin{demo}{Proof}
By the definition of $\Ga_J(P)$, the coefficients $a_j$ have the right differentiability for the proof of \ref{main}(2) to work. 
Definability was used in the proof of \ref{main}(2) only in (IIc) and in claim \ref{main}(4). 
The case (IIc) does not occur, since $E^{(\infty)}(P)=\emptyset$.
In claim \ref{main}(4), the use of definability can be replaced by the following argument:
If a real valued $C^p$ function $f$ vanishes on $t_k \searrow t_0$, then $f^{(q)}(t_0)=0$ for all $0 \le q \le p$.
This follows from a repeated application of Rolle's theorem.
\qed\end{demo}

\begin{lemma} \label{splitdiff}
Let $I \subseteq \R$ be an open interval containing $0$. Let $p, r \in \N_{>0}$.
Suppose that $a_k(t) = t^{kr} a_{(r),k}(t) \in C^{p+nr}(I)$, for $2 \le k \le n$, $a_{(r),2}(0) \ne 0$, and consider 
\[
P_{(r)}(t)(x) = x^n + \sum_{k=2}^n (-1)^k a_{(r),k}(t) x^{n-k}. 
\]
Factorize $P_{(r)} = \prod_{j=1}^l P_{(r),j}$ near $0$, according to the splitting lemma \ref{split}, such that 
\[
P_{(r),j}(t)(x) = x^{n_j} + \sum_{k=1}^{n_j} (-1)^k a_{(r),j,k}(t) x^{n_j-k}, \quad 1 \le j \le l, 
\]
and the $P_{(r),j}$ have mutually distinct roots.
Then, for all $1 \le j \le l$ and $1 \le k \le n_j$, 
$a_{j,k}(t):=t^{kr} a_{(r),j,k}(t)$ belongs to $C^{p+k r}$ near $0$.
\end{lemma}

\begin{demo}{Proof}
By assumption, $t^m a_{(r),k}(t) \in C^{p+(n-k)r+m}$, for all $2 \le k \le n$ and $0 \le m \le kr$.
We assert that 
\begin{equation} \label{eqdiff1}
t^m a_{(r),k}^{(m)}(t) \in C^{p+(n-k)r}, \quad \text{for all $2 \le k \le n$ and $0 \le m \le kr$}, 
\end{equation}
(where $a_{(r),k}^{(m)}$ is understood as distributional derivative).
This follows from
\[
\p_t^m (t^m a_{(r),k}(t)) = \sum_{j=0}^m \binom{m}{j} \frac{m!}{j!} t^j a_{(r),k}^{(j)}(t)
\]
and from induction on $m$.
From \eqref{eqdiff1} we can deduce in a similar way that 
\begin{equation} \label{eqdiff2}
t^q a_{(r),k}^{(q)}(t) \in C^{p}, \quad \text{for all $2 \le k \le n$ and $0 \le q \le nr$}.
\end{equation}

Let $a_{(r)} := (a_{(r),2},\ldots,a_{(r),n})$.
By assumption, there exist $C^\om$ functions $\Ph_{j,k}$ defined in a neighborhood of 
$a_{(r)}(0) \in \R^{n-1}$ such that $a_{(r),j,k} = \Ph_{j,k} \o a_{(r)}$, for all $1 \le j \le l$ and $1 \le k \le n_j$.
Then
\begin{align*}
a_{j,k}^{(k r)}(t) = \sum_{m=0}^{k r} \binom{k r}{m} \frac{(kr)!}{m!}A_{j,k}^m(t),  
\end{align*}
where (by Fa\`a di Bruno, \cite{FaadiBruno1855} for the 1-dimensional version)
\begin{align*}
A_{j,k}^m(t) = \sum_{l\ge 0} \sum_{\substack{\al\in \N_{>0}^l\\ \al_1+\dots+\al_l =m}}
\frac{m!}{l!}d^l \Ph_{j,k}(a_{(r)}(t))\Big( 
\frac{t^{\al_1} a_{(r)}^{(\al_1)}(t)}{\al_1!},\dots,
\frac{t^{\al_l} a_{(r)}^{(\al_l)}(t)}{\al_l!}\Big).
\end{align*}
So, by \eqref{eqdiff2}, we find $a_{j,k}^{(k r)} \in C^p$ and, thus, $a_{j,k} \in C^{p+k r}$.
\qed\end{demo}

\begin{lemma} \label{distinct}
Adopt the setting of lemma \ref{splitdiff}. However, assume that $a_k(t) = t^{kr} a_{(r),k}(t) \in C^{p+kr}(I)$, for $2 \le k \le n$, and that 
all roots of $P_{(r)}(0)$ are distinct.
If $\la_j$ are $C^p$ functions representing the roots of $P_{(r)}$, then $\La_j(t) := t^r \la_j(t)$ are $C^{p+r}$ functions representing the 
roots of $P$.
\end{lemma}

\begin{demo}{Proof}
Instead of \eqref{eqdiff1} we obtain 
\begin{equation} \label{eqdiff3}
t^m a_{(r),k}^{(m)}(t) \in C^{p}, \quad \text{for all $2 \le k \le n$ and $0 \le m \le kr$}.
\end{equation}
The second part of the proof is the same as in \ref{splitdiff}, where now $l=n$ and $n_j=1$ for all $j$. 
In the end we use \eqref{eqdiff3} instead of \eqref{eqdiff2}.
\qed\end{demo}

\subsection{}
Let $P(t)$, $t \in I$, be a curve of monic hyperbolic polynomials of degree $n$ (not necessarily definable). 
Assume $E^{(\infty)}(P)=\emptyset$.
Let $t_0 \in I$ and suppose that the coefficients of $P$ belong to $C^{\Ga_{t_0}(P)+1}$ near $t_0$.
The gradual splitting of, firstly, $P$ near $t_0$ into factors $P_i$ with mutually distinct roots such that all roots of $P_i(t_0)$ coincide, 
then, secondly, of each $(P_i)_{(r_i)}$ (defined in \eqref{Pr}) and so on, determines a well-defined mapping $(P,t_0) \mapsto T(P,t_0)$, where 
$T(P,t_0)$ is a rooted tree in $\cR(n)$ (cf.\ \ref{tree}). 

By the \emph{height} $h(T)$ of a tree $T$ we mean the maximal length (number of edges) of paths connecting the root with a leaf in $T$. 
The \emph{$k$-level} of $T$ is the set of all vertices whose distance (length of the connecting path) from the root is $k$.

\begin{figure}[htp]
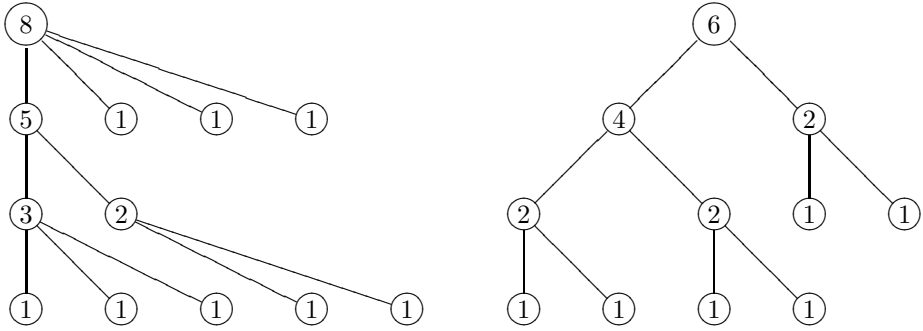

\[
\subfigure{
\xygraph{
[] *+[o]+[F]{8}(
-[d] *+=[o]+[F]{5}(
-[d] *+=[o]+[F]{3}(
-[d] *+=[o]+[F]{1},
-[dr] *+=[o]+[F]{1},
-[drr] *+=[o]+[F]{1}
),
-[dr] *+=[o]+[F]{2}(
-[drr] *+=[o]+[F]{1},
-[drrr] *+=[o]+[F]{1}
),
),
-[dr] *+=[o]+[F]{1},
-[drr] *+=[o]+[F]{1},
-[drrr] *+=[o]+[F]{1},
)
}
}
\hspace{1cm}
\subfigure{
\xygraph{
[] *+[o]+[F]{6}(
-[dl] *+=[o]+[F]{4}(
-[dl] *+=[o]+[F]{2}(
-[d] *+=[o]+[F]{1},
-[dr] *+=[o]+[F]{1}
),
-[dr] *+=[o]+[F]{2}(
-[d] *+=[o]+[F]{1},
-[dr] *+=[o]+[F]{1}
),
),
-[dr] *+=[o]+[F]{2}(
-[d] *+=[o]+[F]{1},
-[dr] *+=[o]+[F]{1}
),
)
}
}
\]
\caption{The first rooted tree is of type (A), the second is not.}
\end{figure}

\begin{theorem} \label{fidiff}
Let $I \subseteq \R$ be an open interval and let $J \subseteq I$ be a compact subinterval. 
Consider a curve of hyperbolic polynomials
\[
P(t)(x) = x^n + \sum_{j=1}^n (-1)^j a_j(t) x^{n-j}, \quad (t \in I)
\]
such that $E^{(\infty)}(P)=\emptyset$.
Assume that the following condition is satisfied for all $t \in J$:
\begin{enumerate}
\item[(A)] For all $k \le h(T(P,t))-2$, the $k$-level of $T(P,t)$ contains at most one vertex with label $\ge 2$.
\end{enumerate}
For each $p \in \N_{>0}$ we have:
\begin{enumerate}
\item[(1)] If the $a_j$ are $C^{p+\Ga_J(P)}$,
then the roots of $P$ can be parameterized by $C^{p+\ga_J(P)}$ functions, globally near $J$.
\end{enumerate}
\end{theorem}

\begin{demo}{Proof}
By \ref{main}(4) and the argument in the proof of \ref{Cp}, it suffices to show the local assertion.
Let $t_0 \in J$ be fixed. 

\begin{claim*}[2]
If the $a_j$ are $C^{p+\Ga_{t_0}(P)}$, then the roots of $P$ can be chosen in $C^{p+\ga_{t_0}(P)}$, locally near $t_0$.
\end{claim*}

Without loss assume that $0 \in J$ and $t_0=0$.
We proceed by induction on $n$. If $n=1$ then $\Ga_0(P)=\ga_0(P)=0$ and we are done.
Suppose $n>1$ and the claim is proved for degrees $\le n-1$.

(I) If $P(0)$ has distinct roots, we have a factorization $P(t)=P_1(t) \cdot P_2(t)$ near $0$, by the splitting lemma \ref{split}.
The coefficients of each factor $P_i$ belong to $C^{p+\Ga_{0}(P)}$. Let $p_i := p+\Ga_{0}(P)-\Ga_{0}(P_i)$. Then $p_i \ge p$, by \eqref{Ga_split}.
Clearly, each $T(P_i,0)$ is of type (A).
By the induction hypothesis, the roots of $P_i$ admit a local parameterization in $C^{p_i+\ga_{0}(P_i)}$. 
By \eqref{ga_split}, $p_i+\ga_{0}(P_i) \ge p+\ga_{0}(P)$, hence claim (2).

(II) If all roots of $P(0)$ coincide, we reduce to the case $a_1=0$. So $a_2(0)=0$.
If $a_2=0$ identically, then all roots are $0$ identically, and claim (2) is satisfied. Suppose that $a_2 \ne 0$.
Since $E^{(\infty)}(P)=\emptyset$ and since $\Ga_0(P)\ge m_0(a_2)$ by definition, we have $m_0(a_2) =2r < \infty$.
Consider $P_{(r)}$ as in \eqref{Pr}.
The coefficients of $P_{(r)}$ are in $C^{p+\Ga_0(P)-nr}$, and $a_{(r),2}(0)\ne 0$. 
Factorize $P_{(r)}(t) = P_{(r),1}(t) \cdots P_{(r),l}(t)$ near $0$ according to the splitting lemma 
\ref{split}. Let $p_{(r),j} := p + \Gamma_0(P_{(r)}) - \Gamma_0(P_{(r),j})$. By the induction hypothesis, 
there are $C^{p_{(r),j} + \gamma_0(P_{(r),j})}$ functions which represent the roots of $P_{(r),j}$ near $0$. 
Let us denote the collection of these functions, for $1 \le j \le l$, by $\la_1,\ldots,\la_n$.  
Then the functions $\La_j(t)=t^r \la_j(t)$ form a local parameterization of the roots of $P$. 
The proof of claim (2) is complete once claim (3) below is shown.

\begin{claim*}[3]
Each $\La_j$ belongs to $C^{p+\ga_0(P)}$.
\end{claim*}

We treat the following cases separately:

(3a) Suppose that $h(T(P,0))\le 2$.

If all $\la_j(0)$ are distinct, then claim (3) follows from \eqref{ga_(r)} and lemma \ref{distinct}.

Otherwise, we can assume (after possibly reordering the $\la_j$) that 
\[
\la_1(0)\!=\!\cdots\!=\! \la_{n_1}(0) \!<\! \la_{n_1+1}(0)\!=\!\cdots\!=\!\la_{n_1+n_2}(0) \!< \!\cdots\! <\! \la_{n-n_l+1}(0)\!=\!\cdots\!=\!\la_n(0).
\]
Set $N(1):=0$ and $N(j):=n_1 +\cdots + n_{j-1}$ for $2 \le j \le l$.
By the splitting lemma \ref{split}, for each $1 \le j \le l$, 
\[
P_{(r),j}(t)(x)=x^{n_j} + \sum_{k=1}^{n_j} (-1)^k a_{(r),j,k}(t) x^{n_j-k}:=\prod_{i=1}^{n_j} (x-\la_{N(j)+i}(t))
\]
has $C^{p+\Ga_0(P_{(r)})}$ coefficients $a_{(r),j,k}$ near $0$. 
By replacing the variable $x$ by $x-a_{(r),j,1}(t)/n_j$, we obtain
\[
\bar P_{(r),j}(t)(x)=x^{n_j} + \sum_{k=2}^{n_j} (-1)^k \bar a_{(r),j,k}(t) x^{n_j-k}=\prod_{i=1}^{n_j} (x-(\la_{N(j)+i}(t)-\tfrac{a_{(r),j,1}(t)}{n_j})),
\]
where the $\bar a_{(r),j,k}$ are still $C^{p+\Ga_0(P_{(r)})}$ near $0$. 
All roots of $\bar P_{(r),j}(0)$ are equal to $0$. 
As above we may conclude that
there is a $q_j \in \N_{>0}$ such that $\bar a_{(r),j,k}(t) = t^{k q_j} \bar a_{(r,q_j),j,k}(t)$, for $2 \le k \le n_j$, $\bar a_{(r,q_j),j,2}(0) \ne 0$, and 
\[
\bar P_{(r,q_j),j}(t)(x):=x^{n_j} + \sum_{k=2}^{n_j} (-1)^k  \bar a_{(r,q_j),j,k}(t) x^{n_j-k}
\]
has $C^{p_{(r),j}+\Ga_0(\bar P_{(r,q_j),j})}$ coefficients $\bar a_{(r,q_j),j,k}(t)$ and $C^{p_{(r),j}+\ga_0(\bar P_{(r,q_j),j})}$ roots $\mu_{j,i}$.
Then
\begin{equation} \label{mula}
t^{q_j} \mu_{j,i}(t) = \la_{N(j)+i}(t)-\tfrac{a_{(r),j,1}(t)}{n_j}, \quad \text{for } 1 \le i \le n_j.
\end{equation}
Thus,
\begin{equation} \label{muLa}
\La_{N(j)+i}(t) = t^{r+q_j} \mu_{j,i}(t) + t^r \tfrac{a_{(r),j,1}(t)}{n_j}, \quad \text{for } 1 \le i \le n_j.
\end{equation}

By lemma \ref{splitdiff},
\[
t^{kr} a_{(r),j,k}(t) \in C^{p+\Ga_0(P_{(r)})+kr}, \quad \text{for all $1 \le j \le l$ and $1 \le k \le n_j$}.
\]
In particular, $t^{r} a_{(r),j,1}(t) \in C^{p+\Ga_0(P_{(r)})+r}$. So, in order to show claim (3), it remains to prove that the first summand 
on the right-hand side of \eqref{muLa} belongs to $C^{p+\ga_0(P)}$. 

The mapping $(a_{(r),j,1},\ldots,a_{(r),j,n_j}) \mapsto (\bar a_{(r),j,2},\ldots,\bar a_{(r),j,n_j})$ is polynomial.
Thus, there exist $C^\om$ functions $\bar \Ph_{j,k}$ defined in a neighborhood of 
$a_{(r)}(0) \in \R^{n-1}$ such that $\bar a_{(r),j,k} = \bar \Ph_{j,k} \o a_{(r)}$, for all $1 \le j \le l$ and $2 \le k \le n_j$.
Hence, by (the proof of) lemma \ref{splitdiff}, we also obtain
\[
t^{kr} \bar a_{(r),j,k}(t) \in C^{p+\Ga_0(P_{(r)})+kr}, \quad \text{for all $1 \le j \le l$ and $2 \le k \le n_j$},
\]
and thus (by \eqref{Ga_(r)})
\begin{align*}
t^{k(r+q_j)} \bar a_{(r,q_j),j,k}(t) \in C^{p+\Ga_0(P_{(r)})+kr} &\subseteq C^{p_{(r),j}+\Ga_0(\bar P_{(r,q_j),j})+k(r+q_j)}, \\
&\quad \text{for all $1 \le j \le l$ and $2 \le k \le n_j$}.
\end{align*}

By the assumption $h(T(P,0)) \le 2$, all $\mu_{j,i}(0)$ are distinct. Then claim (3) follows from \eqref{ga_split}, \eqref{ga_(r)}, and lemma \ref{distinct}.

(3b) Suppose that $h(T(P,0)) > 2$. Let us use the notation of (3a). Since $T(P,0)$ is of type (A), we may assume $n_2=n_3=\cdots=n_l=1$,
and the roots $\la_j$, for $n_1+1 \le j \le n$, belong to $C^{p+\Ga_0(P_{(r)})}$.  

Consider the Newton polynomials $s_{(r),k} = \sum_{j=1}^n \la_j^k$ and $\bar s_{(r,q_1),1,k} = \sum_{j=1}^{n_1} \mu_{1,j}^k$,
associated to $P_{(r)}$ and $\bar P_{(r,q_1),1}$, respectively.
(In the following argument it is convenient to work with the Newton polynomials of the roots instead of the elementary symmetric functions (coefficients).
They are related to each other by the polynomial diffeomorphism defined in \eqref{rec}.)
Note that $s_{(r),1}=\bar s_{(r,q_1),1,1}=0$ and $\bar s_{(r,q_1),1,0} = n_1$.
We have, by \eqref{mula},
\begin{gather}
0=s_{(r),1} = a_{(r),1,1}(t) + \sum_{i=n_1+1}^n \la_i(t), \label{bla}\\
s_{(r),k}(t) = \sum_{i=0}^k \binom{k}{i} t^{iq_1} \bar s_{(r,q_1),1,i}(t) \big(\tfrac{a_{(r),1,1}(t)}{n_1}\big)^{k-i} +\sum_{i=n_1+1}^n \la_i(t)^k , 
\quad 2 \le k \le n_1. \label{su}  
\end{gather}
By lemma \ref{splitdiff} and \eqref{Ga_(r)}, $\La_i(t) = t^r \la_i(t) \in C^{p+\Ga_0(P)}$, for $n_1+1 \le i \le n$. 
Thus, by \eqref{bla}, $t^r a_{(r),1,1}(t) \in C^{p+\Ga_0(P)}$.
By \eqref{rec}, we have $t^{kr} s_{(r),k}(t) \in C^{p+\Ga_0(P)}$, for $2 \le k \le n$ 
(since the same is true when the $s_{(r),k}$ are replaced by the $a_{(r),k}$). 
Hence, \eqref{su} implies inductively that $t^{i(r+q_1)} \bar s_{(r,q_1),1,i}(t) \in C^{p+\Ga_0(P)}$, for $2 \le i \le n_1$, and equivalently, 
\[
t^{i(r+q_1)} \bar a_{(r,q_1),1,i}(t) \in C^{p+\Ga_0(P)},\quad \text{for } 2 \le i \le n_1.
\]
Let us repeat this procedure with
\[
\tilde P(t)(x):=x^{n_1}+\sum_{j=2}^{n_1} (-1)^j t^{j(r+q_1)} \bar a_{(r,q_1),1,j}(t) x^{n_1-j} = \prod_{i=1}^{n_1} (x-t^{r+q_1}\mu_{1,i}(t))
\]
instead of $P$. Evidently, $T(\tilde P,0)$ is of type (A).
After finitely many steps the situation is reduced to case (3a).
This completes the proof of claim (3).

\begin{claim*}[4]
If the $a_j$ are $C^{p+\Ga_J(P)}$, then the roots of $P$ can be chosen in $C^{p+\ga_J(P)}$, locally near $t_0$.
\end{claim*}

By claim (2), the roots of $P$ can be chosen in $C^{p+\Ga_J(P)-\Ga_{t_0}(P)+\ga_{t_0}(P)}$, locally near $t_0$.
By definition, $p+\Ga_J(P)-\Ga_{t_0}(P)+\ga_{t_0}(P) \ge p+\ga_J(P)$.
\qed\end{demo}

\begin{remark}
We do not know whether or not theorem \ref{fidiff} holds, if $T(P,t)$ is not of type (A).
Note that each $T \in \bigcup_{n=1}^4 \cR(n)$ is automatically of type (A).
Thus, theorem \ref{fidiff} is true for all $P$ with degree at most $4$.
\end{remark}

\section{Definable version of Bronshtein's theorem}

\begin{theorem} \label{B}
Let $I \subseteq \R$ be an open interval. Consider a curve of hyperbolic polynomials
\[
P(t)(x) = x^n + \sum_{j=1}^n (-1)^j a_j(t) x^{n-j}
\]
with definable $C^n$ coefficients $a_j$.
Then the roots of $P$ can be parameterized by definable $C^1$ functions, globally.
\end{theorem}

If `definable' is omitted in the formulation of theorem \ref{B}, then we obtain Bronshtein's theorem \cite{Bronshtein79} (see also \cite{Wakabayashi86}). 
Actually we obtain the refinement of Bronshtein's theorem due to \cite{ColombiniOrruPernazza08}.
The proof of Bronshtein's theorem is very delicate and only poorly understood. In the definable case it becomes remarkably simple. 

\begin{demo}{Proof}
By \ref{main}(4), it suffices to show the local statement.
We follow the proof of theorem \ref{main}(3) and indicate the necessary modifications.
Let us begin the induction on $n$ with the case $n=1$, which is trivial.
(I) and (II) can be adopted with obvious minor changes. So assume that $a_1=0$ identically and $a_2(0)=0$.
Since $0 \le \tilde \De_2 \o P= -2n a_2$, we have $m_0(a_2) \ge 2$.
By the multiplicity lemma \ref{mult} (for $r=1$), $m_0(a_k)\ge k$ for $2 \le k \le n$, and $P_{(1)}$ (defined in \eqref{Pr}) 
is a continuous curve of hyperbolic polynomials. Let $\mu_j$ be a continuous 
parameterization of the roots of $P_{(1)}$ near $0$.  
Then the functions $\la_j(t):=t \mu_j(t)$ form a definable continuous parameterization of 
the roots of $P$ near $0$ such that $m_0(\la_j)\ge 1$ for each $j$.
By lemma \ref{C^p-def}(1), each $\la_j$ is $C^1$ near $0$. 
\qed\end{demo}

\begin{examples}
(1) The function $f(t) = t^2 |t|$ is in $C^{2,1}$ (but not three times differentiable). 
The square roots of $f$ may be chosen $C^1$ but not $C^{1,1}$.

(2) Let $g(t) = 1/3$ for $t \ge 0$ and $g(t) = 0$ otherwise. 
Consider the following $C^{2,1}$ curve of cubic polynomials
(cf.\ \cite[Example 4.6]{ColombiniOrruPernazza08}):
\[
P(t)(x) = x^3- t^3 g(t) x^2 + (2 t^3 g(t)-t^2) x - t^3 g(t).
\]
Its discriminant is $\tilde \De_3(P(t)) = t^6 (1+o(1))$ if $t\ge 0$ 
and $\tilde \De_3(P(t)) =4 t^6$ if $t < 0$.
Thus, for small $t$, $P(t)$ is hyperbolic.
The roots of $P$ cannot be chosen differentiable at $0$:
Note that $0$ is a triple root of $P(0)$. Consider, for $t\ne 0$, 
\[
Q(t)(y) = t^{-3}P(t)(ty) = y^3- t^2 g(t) y^2 + (2 t g(t)- 1)y-g(t).
\]
Then $\lim_{t \searrow 0} Q(t)(y) = y^3-y-1/3$ and $\lim_{t \nearrow 0} Q(t)(y) = y^3-y$.
Thus, the roots of $P$ cannot be differentiable at $0$.
\end{examples}

\section{Complex polynomials}

\subsection{} 
In this section 
let $I \subseteq \R$ be an open interval
and consider a curve of complex polynomials
\[
P(t)(x) = x^n + \sum_{j=1}^n (-1)^j a_j(t) x^{n-j},
\]
i.e., each coefficient $a_j : I \to \C$ is a continuous complex valued function.
Then the roots of $P$ admit a continuous parameterization (e.g.\ \cite[II 5.2]{Kato76}).

A complex valued function $f : I \to \C$ is called \emph{definable} if $(\on{Re} f, \on{Im} f) : I \to \R^2$ is definable.
We will assume that the coefficients $a_j$ of $P$ are definable.

The set $E^{(\infty)}(P)$ can be defined and has the same properties as in the hyperbolic case (cf.\ \ref{Einfty}).

\begin{lemma} \label{Cdef}
If the coefficients $a_j$ of $P$ are definable, then every continuous parameterization $\la_j$ of the roots of $P$ is definable.
\end{lemma}

\begin{demo}{Proof}
The real and imaginary parts $\on{Re} \la_j$, $\on{Im} \la_j$,  $1 \le j \le n$, 
parameterize the solutions of the $2n$ algebraic equations with definable coefficients
$\on{Re} P(t)(\la_j(t))=0$, $\on{Im} P(t)(\la_j(t))=0$, $1 \le j \le n$.
The family of continuous parameterizations of the solutions of these equations is finite.
\qed\end{demo}

\begin{theorem} \label{complex}
Let $I \subseteq \R$ be an open interval. Consider a curve of polynomials
\[
P(t)(x) = x^n + \sum_{j=1}^n (-1)^j a_j(t) x^{n-j}, 
\]
with definable $C^\infty$ coefficients $a_j$.
Then, for each $t_0 \in I$, there is an $N \in \N_{>0}$ such that
$t \mapsto P(t_0 \pm (t-t_0)^N)$ admits definable $C^\infty$ parameterizations of its roots, locally near $t_0$.
\end{theorem}

\begin{demo}{Proof}
 Since the coefficients of $t \mapsto P(t_0 \pm (t-t_0)^N)$ are definable, we need not care about the definability of its roots, by lemma \ref{Cdef}.
Without loss assume that $0 \in I$ and $t_0=0$. 
We proceed by induction on $n$. The case $n=1$ is trivial. 

(I) If $P(0)$ has distinct roots, we are done, by the splitting lemma \ref{split} and the induction hypothesis.
(Here we use that, if $t \mapsto P_i(\pm t^{N_i})$, $i=1,2$, admit $C^\infty$ roots then so does $t \mapsto P_1(\pm t^{N_1 N_2}) P_2(\pm  t^{N_1 N_2})$.)

(II) If all roots of $P(0)$ coincide, we reduce to the case $a_1=0$. Then all roots of $P(0)$ are equal to $0$. 

(IIa) If $m_0(a_k) < \infty$ for some $2 \le k \le n$, there exist $N,r \in \N_{>0}$ such that $(t \mapsto P(\pm t^N))_{(r)}$ 
(the reduced curve of polynomials defined in \eqref{Pr} associated to $t \mapsto P(\pm t^N)$)
has distinct roots at $t=0$ (see \cite{RainerAC}). By the splitting lemma \ref{split} and the induction hypothesis, we are done. 

(IIb) If all $a_k=0$ identically, then all roots of $P$ are identically $0$.

(IIc) If $m_0(a_k) = \infty$ for all $2 \le k \le n$, then for any continuous choice $\la_j$ of the roots of $P$ we find $m_0(\la_j)=\infty$ (for all $j$).
For: Let $\la(t)$ be any continuous root of $P(t)$ and $r \in \N_{>0}$.
Then, for $t\ne 0$, $\mu(t) = t^{-r} \la(t)$ is a root of $P_{(r)}(t)$ (defined in \eqref{Pr}), hence bounded in $t$.
So $\la(t) = t^{r-1} \cdot t \mu(t)$, and $t \mapsto t \mu(t)$ is continuous.
Thus $m_0(\la)=\infty$, since $r$ was arbitrary.
By lemma \ref{C^p-def}(1) (applied to $\on{Re} \la_j$ and $\on{Im} \la_j$), for each $p$, 
there is a neighborhood $I_p$ of $0$ such that each $\la_j$ is $C^p$ on $I_p$.
Since the coefficients $a_j$ (and hence the $\tilde \De_k \o P$) are definable,
for small $t \ne 0$ the multiplicity of the $\la_j(t)$ is constant.
So all $\la_j$ are $C^\infty$ off $0$ (by the splitting lemma \ref{split}) and hence also near $0$.
\qed\end{demo}

\subsection{}
In \cite{RainerAC} we have deduce from the analog of theorem \ref{complex} that any continuous parameterizations of the roots of 
a $C^\infty$ curve $P$ of complex polynomials with $E^{\infty}(P)=\emptyset$ is locally actually absolutely continuous 
(not better!, see \ref{bestac} below). 
The optimal conditions for absolutely continuous roots are unknown.
 
However, in the definable case we have the following best possible result:

\begin{theorem}
Any continuous choice of the roots of a curve of monic complex polynomials with definable continuous coefficients is locally absolutely continuous. 
\end{theorem}

\begin{demo}{Proof}
This follows from lemma \ref{Cdef} and lemma \ref{defAC} below.
\qed\end{demo}

\begin{lemma} \label{defAC}
Let $I \subseteq \R$ be an interval.
A definable continuous function $f : I \to \C$ is locally absolutely continuous.
\end{lemma}

\begin{demo}{Proof}
We show that a continuous definable function $f : I \to \R$, where $I \subseteq \R$ is a compact interval, is absolutely continuous. 
By the Monotonicity theorem \cite{vandenDries98}, $f$ is $C^1$ on the complement of finitely many points $J = I \setminus \{a_1,\ldots,a_n\}$.
Let $J_i$ be some connected component of $J$. By definability, we can partition $J_i$ into finitely many subintervals $J_{ij}$ on
each of which either $f' > 0$ or $f' \le 0$. 
Then it is easy to see that $f'|_{J_{ij}}$ belongs to $L^1$ for every $J_{ij}$, thus $f'|_{J_i}$ belongs to $L^1$ (here we use that $f$ is continuous).
Let $[a,b]:=\overline J_i$ denote the closure of $J_i$.
Then we have $f(x)=f(a) + \int_a^x f'(t) dt$ for $x \in [a,b]$.
So $f|_{\overline J_i}$ is absolutely continuous. Since $J_i$ was arbitrary, the proof is complete.
\qed\end{demo}

\begin{examples} \label{bestac}
Absolute continuity is the best we can hope for:
In general the roots cannot be chosen with first derivative in $L^p_{\on{loc}}$
for any $1 < p \le \infty$.
This is demonstrated by
\[ 
P(t)(z) = z^n- t, \quad  t \in \R,
\]
for $1 < p < \infty$
if $n \ge \frac{p}{p-1}$
and for $p=\infty$ if $n \ge 2$. 
\end{examples}


\begin{thebibliography}{10}

\bibitem{AKLM98}
D.~Alekseevsky, A.~Kriegl, M.~Losik, and P.~W. Michor, \emph{Choosing roots of
  polynomials smoothly}, Israel J. Math. \textbf{105} (1998), 203--233.

\bibitem{BBCP06}
J.-M. Bony, F.~Broglia, F.~Colombini, and L.~Pernazza, \emph{Nonnegative
  functions as squares or sums of squares}, J. Funct. Anal. \textbf{232}
  (2006), no.~1, 137--147. 

\bibitem{Bronshtein79}
M.~D. Bronshtein, \emph{Smoothness of roots of polynomials depending on
  parameters}, Sibirsk. Mat. Zh. \textbf{20} (1979), no.~3, 493--501, 690,
  English transl. in Siberian Math. J. \textbf{20} (1980), 347--352. See also
  MR0537355 (82c:26018). 

\bibitem{ColombiniOrruPernazza08}
F.~Colombini, N.~Orr\'u, and L.~Pernazza, \emph{On the regularity of the roots
  of hyperbolic polynomials}, Preprint, 2008.

\bibitem{FaadiBruno1855}
C.F. Fa\`a~di Bruno, \emph{Note {s}ur {u}ne {n}ouvelle {f}ormule {d}u {c}alcul
  {d}iff\'erentielle}, Quart. J. Math. \textbf{1} (1855), 359--360.

\bibitem{Glaeser63R}
G.~Glaeser, \emph{Racine carr\'ee d'une fonction diff\'erentiable}, Ann. Inst.
  Fourier (Grenoble) \textbf{13} (1963), no.~fasc. 2, 203--210. 

\bibitem{Kato76}
T.~Kato, \emph{Perturbation theory for linear operators}, second ed.,
  Grundlehren der Mathematischen Wissenschaften, vol. 132, Springer-Verlag,
  Berlin, 1976. 

\bibitem{KLM04}
A.~Kriegl, M.~Losik, and P.~W. Michor, \emph{Choosing roots of polynomials
  smoothly. {II}}, Israel J. Math. \textbf{139} (2004), 183--188. 

\bibitem{Procesi78}
C.~Procesi, \emph{Positive symmetric functions}, Adv. in Math. \textbf{29}
  (1978), no.~2, 219--225. 

\bibitem{RainerAC}
A.~Rainer, \emph{Perturbation of complex polynomials and normal operators},
  Math. Nach. \textbf{282} (2009), No.~12, 1623--1636. 

\bibitem{Rellich37I}
F.~Rellich, \emph{St\"orungstheorie der {S}pektralzerlegung}, Math. Ann.
  \textbf{113} (1937), no.~1, 600--619. 

\bibitem{vandenDries98}
L.~van~den Dries, \emph{Tame topology and o-minimal structures}, London
  Mathematical Society Lecture Note Series, vol. 248, Cambridge University
  Press, Cambridge, 1998. 

\bibitem{Wakabayashi86}
S.~Wakabayashi, \emph{Remarks on hyperbolic polynomials}, Tsukuba J. Math.
  \textbf{10} (1986), no.~1, 17--28. 

\end{thebibliography}

\def\cprime{$'$}
\providecommand{\bysame}{\leavevmode\hbox to3em{\hrulefill}\thinspace}
\providecommand{\MR}{\relax\ifhmode\unskip\space\fi MR }
\providecommand{\MRhref}[2]{%
  \href{http://www.ams.org/mathscinet-getitem?mr=#1}{#2}
}
\providecommand{\href}[2]{#2}

\end{document}